\documentclass{conm-p-l}

\usepackage{tikz}
\usetikzlibrary{shapes.symbols,arrows}

\copyrightinfo{2014}{Thomas Nowak and Bernadette Charron-Bost}

\newtheorem{theorem}{Theorem}[section]
\newtheorem{lemma}[theorem]{Lemma}

\theoremstyle{definition}

\theoremstyle{remark}

\numberwithin{equation}{section}

\DeclareMathOperator{\ind}{ind}
\DeclareMathOperator{\Red}{Red}
\newcommand{\IN}{{\mathbb{N}}}
\newcommand{\IR}{{\mathbb{R}}}
\newcommand{\IRmax}{\IR_{\max}}
\newcommand{\bzero}{-\infty}
\newcommand{\bone}{0}

\begin{document}

\title{An Overview of Transience Bounds in Max-Plus Algebra}


\author{Thomas Nowak}
\address{Laboratoire d'Informatique, \'Ecole polytechnique, 91128 Palaiseau,
France}
\curraddr{}
\email{nowak@lix.polytechnique.fr}
\thanks{}

\author{Bernadette Charron-Bost}
\address{CNRS, Laboratoire d'Informatique, \'Ecole polytechnique, 91128 Palaiseau,
France}
\curraddr{}
\email{charron@lix.polytechnique.fr}
\thanks{}


\date{}

\begin{abstract}
We survey and discuss upper bounds on the length of the transient phase of
max-plus linear systems and sequences of max-plus matrix powers.
In particular, we explain how to extend a result by Nachtigall to yield a new
approach for proving such bounds
and we state an asymptotic tightness result by using an example given by
Hartmann and Arguelles.
\end{abstract}

\maketitle


\section{Introduction}

Max-plus linear algebra is used to describe
production and transportation systems, and several distributed algorithms,
due to the occurrence of a temporal maximum operation when events are
synchronized.
These systems are described by the repeated application of a fixed matrix to
an initial vector.
A body of research spawned examining the behavior and parameters
of interest of max-plus linear systems.
If the system matrix is irreducible, i.e., if the digraph described by 
it is strongly connected,
 one observes a periodic behavior after an initial
transient phase whose length we refer to as the {\em transient\/}:
If~$x(k)$ denotes the system vector after~$k$ applications, then
$x_i(k+p) = x_i(k) + c$ with a constant~$c$ independent
of index~$i$ for all~$k$ large enough, i.e., greater or equal to the transient.
This was first shown by Cohen et al.~\cite{CDQV83}.
In fact, they showed that the entries in the sequence of max-plus
powers~$A^{\otimes k}$ of
every irreducible matrix~$A$ are eventually periodic in the same sense, i.e.,
$A_{i,j}^{\otimes (k+p)} = A^{\otimes k}_{i,j} + c$ with some~$c$ independent
of the index~$(i,j)$ for all~$k$ large enough.
This obviously implies the result for systems.
In the same vein, the transient of a system is always upper bounded
by the transient of the sequence of powers of its system matrix.

It is the purpose of this paper to survey the existing upper bounds on the 
transient of max-plus linear systems and matrices.
Preceding the first general transience bound by Hartmann and
Arguelles~\cite{HA99}, a number of upper
bounds on the transient of certain max-plus linear systems
in computer science  were established
(e.g., \cite{ER97,MR92,CBFWW11,Chr84}).
All these results are, as far as we are aware of them, covered by the bounds we chose to
present here.

Hartmann and Arguelles~\cite{HA99} proved, as a corollary of their upper bound
on the transient, that the computation of the exact value of the transient of a
system or a matrix can be done in time polynomial in the size of a list
representation.
Their algorithm first calculates an upper bound on the transient and then
identifies the transient by doing a binary search.
Consequently, better upper bounds improve the running time of this algorithm.
Also, bounds
involving certain parameters of the systems or the matrix allow to design for a
small transient.

The paper is structured as follows:
Section~\ref{sec:prelim} defines the basic notions used in the paper and
presents some preliminary results used later.
Section~\ref{sec:boolean} discusses selected upper bounds on the transient for
the special case of Boolean matrices.
In Section~\ref{sec:nachtigall}, we present the decomposition of the sequence
of matrix powers as a maximum of matrices with bounded transients established by
Nachtigall~\cite{Nac97}.
We also explain how his arguments can be extended and completed to show a bound on
the transient of the sequence of powers.
This extension was presented by the authors at the workshop
TROPICAL-12~\cite{CBN12}. 
Section~\ref{sec:ha} gives the first general transience bound proved by
Hartmann and Arguelles~\cite{HA99}.
It also includes the generalization of one of their examples showing a form of
asymptotic tightness of their bound, and also of other bounds.
In Section~\ref{sec:bg}, we present a bound for a special class of max-plus
matrices that was proved with a technique used by Akian et al.~\cite{AGW05} and
Bouillard and Gaujal~\cite{BG01}.
Section~\ref{sec:syk} gives the bound by Soto y Koelemeijer~\cite{SyK03} for the special
case of max-plus matrices whose all entries are finite.
Section~\ref{sec:cbfn} presents the bounds by Charron-Bost et al.~\cite{CBFN12}
that use
two graph parameters (girth and cyclicity).
In Section~\ref{sec:sys:mat}, we present the method by Charron-Bost et al.\ to
transform transience bounds for systems into transience bounds for matrices.
With Section~\ref{sec:conclusion}, we conclude the paper by comparing the
bounds from both a quantitative and a qualitative viewpoint.

\section{Definitions and Preliminaries}\label{sec:prelim}

In max-plus algebra, one endows the set~$\IRmax = \IR \cup \{-\infty\}$ with
the addition $a\oplus b = \max\{a,b\}$ and the multiplication $a\otimes b =
a+b$.
The identity with respect to~$\oplus$ is $\bzero$ and~$\bone$ with respect 
to~$\otimes$.

A max-plus matrix is a matrix with entries in~$\IRmax$.
The max-plus product of two matrices~$A$ and~$B$ of compatible size is defined
in the usual way by setting $(A\otimes B)_{i,j} = \bigoplus_h A_{i,h}\otimes B_{h,j}$.
We write~$A^{\otimes k}$
for the $k$th max-plus power of a square matrix~$A$.

A digraph is a pair $G=(V,E)$ of a nonempty set~$V$ of nodes and a set
$E\subseteq V\times V$ of edges.
A walk in the digraph is a {\em path\/} if every node occurs only once.
A closed walk is a {\em cycle\/} if only the start and end node occurs twice.
We write~$\ell(W)$ for the length of walk~$W$.

The length of the shortest cycle in a digraph~$G$ is called the {\em girth\/}
of~$G$.
If a digraph is strongly connected, the greatest common divisor of its cycle
lengths is called its {\em cyclicity}.
The cyclicity of a (possibly not strongly connected) digraph is the least common
multiple of the cyclicities of its strongly connected components.

To every $n\times n$ max-plus matrix~$A$ corresponds a digraph~$G(A)$ with node
set $V=\{1,2,\dots,n\}$ containing an edge $(i,j)$ if and only if
$A_{i,j}\neq\bzero$.
We refer to~$A_{i,j}$ as the weight of edge $(i,j)$.
Matrix~$A$ is {\em irreducible\/} if~$G(A)$ is strongly connected.
If~$W$ is a walk in~$G(A)$, we define its weight~$A(W)$ as the sum of the
weights of its edges.
The entry~$A_{i,j}^{\otimes k}$ is the maximum weight of walks from~$i$ to~$j$
of length~$k$.\footnote{We follow the convention that $\max\emptyset = -\infty$.}
If~$v$ is a max-plus column vector of size~$n$, then
the entry $\big(A^{\otimes k} \otimes v\big)_i$ is the maximum of the values
$A(W)+v_j$ where the maximum is formed over all nodes~$j$ and all walks~$W$
from~$i$ to~$j$ of length~$k$.

Denote by~$\lambda(A)$ the maximum mean weight $A(\Gamma)/\ell(\Gamma)$ of
cycles in~$G(A)$.
We call {\em critical\/} every cycle with maximum mean weight.
The subgraph of~$G(A)$ induced by edges on critical cycles is called the {\em
critical subgraph}.

Let~$X$ be a fixed index set.
A sequence $a:\IN\to\IRmax^X$ is called {\em eventually periodic\/} with
period~$p\geq1$
and ratio~$\varrho\in\IR$ if there exists a~$K\in\IN$ such that
$a(k+p)=a(k)+p\cdot\varrho$ for all $k\geq K$, where the addition is to be
understood component-wise.
The smallest such~$K$ is the {\em transient\/} of the
sequence.
The ratio is unique if not all components of~$a(k)$ are eventually
constant~$\bzero$.
For every eventual periodicity of sequence~$a(k)$ with ratio~$\varrho$,
 the sequence $a(k) - k\cdot\varrho$ is eventually periodic with ratio~$0$.

The following lemma shows that the transient is independent of the considered
period.
For a proof see, for instance, \cite[Lemma 11]{HA99}.

\begin{lemma}\label{lem:transient:indep}
Let~$a(k)$ be a sequence, $p,q$ positive integers, $K_p,K_q$ nonnegative
integers, and~$\varrho$ a real number.
If $a(k+p)=a(k)+p\cdot \varrho$ for all $k\geq K_p$ and $a(k+q)=a(k)+q\cdot
\varrho$ for all $k\geq K_q$,
then $a(k+\gcd(p,q))=a(k) + \gcd(p,q)\cdot \varrho$ for all $k\geq \max\{K_p,K_q\}$.
\end{lemma}

Cohen et al.~\cite{CDQV83} showed that the sequence of powers of an irreducible
max-plus matrix, and hence of all systems with irreducible matrix, are
eventually periodic.
Denote by~$\gamma_c(A)$ the cyclicity of the critical subgraph of~$G(A)$.

\begin{theorem}[Cohen et al., 1983]\label{thm:cyc}
The sequence of powers~$A^{\otimes k}$ of an irreducible square max-plus
matrix~$A$ is eventually periodic with ratio~$\lambda(A)$ and
period~$\gamma_c(A)$.
\end{theorem}

Theorem~\ref{thm:cyc} is based on the fact that maximum weight walk eventually
include in the majority critical cycles.
To give an explicit upper bound on when
they visit at least one critical cycle, several authors defined what they considered to be the ``second most
significant'' cycle mean.
This can be done in a number of ways, depending on the specific proof technique
used:
One possibility, used by the authors' extension of Nachtigall's decomposition
and by Soto y Koelemeijer, is to consider the second largest cycle mean~$\lambda_2(A)$.
Bouillard and Gaujal and Charron-Bost et al.\ both considered the largest cycle mean disjoint to all critical
cycles, which we denote by~$\lambda_{nc}(A)$.
Hartmann and Arguelles used a third parameter, $\lambda_0(A)$,  which is defined in terms of the max-balancing
\cite{SS91} of~$G(A)$.
We do not formally define the three parameters, but give their relative
ordering, also with respect to~$\lambda(A)$:

\begin{equation}\label{eq:lambdas}
\lambda(A) > \lambda_2(A) \geq \lambda_{nc}(A) \geq \lambda_0(A)
\end{equation}

We denote by~$\lVert A\rVert$ the difference between the greatest and smallest
finite entry in matrix~$A$.

\section{The Boolean Case: Index of Convergence}\label{sec:boolean}

A Boolean matrix is a max-plus matrix whose entries are either~$\bzero$
or~$\bone$ and it corresponds to a digraph.
The behavior of the sequence of powers of Boolean matrices, or equivalently the
set of possible walk lengths between nodes in a digraph, has been extensively studied (see, e.g.,
\cite{BR91} or~\cite{LS93} for an overview).
If the digraph is strongly connected, Theorem~\ref{thm:cyc} shows that every
such sequence is eventually periodic; for the critical subgraph is equal to the
whole digraph.
Its transient is commonly referred to as the {\em index of convergence\/} of the matrix
resp.\ the digraph (also sometimes the {\em exponent\/} if the cyclicity is
equal to~$1$).
Clearly, the case of Boolean matrices is an important special case for the
study of transients in max-plus algebra.

The first bound on the index of convergence was given by Wielandt~\cite{Wie50} 
 for the case of {\em primitive\/} digraphs, i.e.,
digraphs whose cyclicity is equal to~$1$.
He also gave a class of examples showing his bound is tight.

\begin{theorem}[Wielandt, 1950]\label{thm:wielandt}
The index of convergence of a strongly connected primitive digraph with~$n$ nodes is at most
$(n-1)^2+1$.
Furthermore, for every~$n\geq2$ there exists a strongly connected primitive digraph with~$n$ nodes whose index of
convergence is equal to $(n-1)^2+1$.
\end{theorem}

The bound of $(n-1)^2+1$ was refined independently by Dulmage and Mendelsohn~\cite{DM64}
and by Denardo~\cite{Den77} in terms of the digraph's girth~$g$.
They arrived at the same bound, which is in the order of $O(g\cdot n)$.
This suggests that the lower the
girth, the lower the index of convergence.

Later, Schwarz~\cite{Sch70} extended Theorem~\ref{thm:wielandt} to
non-primitive digraphs.\footnote{Shao and Li~\cite{SL87} gave an alternative
proof.} 
Interestingly, he showed that the bound of
$(n-1)^2+1$ remains true and that even a lower upper bound holds, which in the
order
of~$O(n^2/\gamma)$ where~$\gamma$ denotes the cyclicity.
This suggests that the higher the cyclicity, the lower the index of
convergence.

Because the girth of a strongly connected digraph is always greater or equal to
the cyclicity, the two results suggest a necessary trade-off between the two parameters
for attaining a small index of convergence.
For instance,
the two parameters need to be equal for attaining the minimal index of convergence of~$0$.

Kim~\cite{Kim79} showed a new upper bound, which generalizes both the bounds of
Dulmage, Mendelsohn, and Denardo, and Schwarz:

\begin{theorem}[Kim, 1979]
The index of convergence of a strongly connected digraph with~$n$ nodes, girth~$g$, and
cyclicity~$\gamma$ is at most
\[
n + g\cdot \left( \left\lfloor \frac{n}{\gamma} \right\rfloor - 2  \right)
\enspace.
\]
\end{theorem}

\section{Nachtigall Decomposition}\label{sec:nachtigall}

A significant step in the direction of a transience bound for non-Boolean
matrices
was done by Nachtigall~\cite{Nac97}.
While he did not prove a bound on the transient, he showed that 
the sequence of matrix powers can be written as a maximum of eventually periodic
sequences with bounded transients.
Such a decomposition in the form of a maximum, by itself, does not yield a bound on the
transient of the original sequence; it does not even imply that it is
eventually periodic.
As a matter of fact, Nachtigall shows the existence
of such a decomposition not only for irreducible matrices, but for general
square max-plus matrices, for which the sequence of powers is not necessarily
eventually periodic.

The authors~\cite{CBN12} have observed that the specific structure of the
Nachtigall decomposition in the case of irreducible matrices allows 
to deduce a transience bound.
We discuss this after the description of the decomposition.
Nachtigall's decomposition was also studied by Moln\'arov\'a~\cite{Mol03} and
Sergeev and
Schneider~\cite{SS12}.

\begin{theorem}[Nachtigall, 1997]\label{thm:nachtigall}
Let~$A$ be an $n\times n$ max-plus matrix.
Then there exist eventually periodic matrix sequences
$A_1(k),A_2(k),\dots,A_n(k)$ with transients at most~$3n^2$ such that
for all $k\geq0$:
\[
A^{\otimes k}  = A_1(k) \oplus A_2(k) \oplus \cdots \oplus A_n(k)
\]
\end{theorem}

Nachtigall proved Theorem~\ref{thm:nachtigall} by recursively picking a
cycle~$\Gamma$ with maximal ratio~$A(\Gamma)/\ell(\Gamma)$ and by partitioning the sets of walks
in~$G(A)$ into the sets of walks that do and do
not visit cycle~$\Gamma$.
Walks that do not visit~$\Gamma$ are walks in the subgraph of~$G(A)$ that has all
edges incident to~$\Gamma$ removed.
This subgraph is the digraph of the matrix obtained from~$A$ by setting
to~$\bzero$ all rows and columns corresponding to nodes in~$\Gamma$; its effective
size is strictly smaller than the size of~$A$, which enables a recursive
descent.
If no cycle exists in~$G(A)$ at all, then the transient of~$A$ is at most~$n$ since in
this case, $A_{i,j}^{\otimes k}=\bzero$ for all~$i,j$ and all $k\geq n$.

One can see that the sequence of the maximum weights of walks of length~$k$
from a node~$i$ to a node~$j$ that do visit cycle~$\Gamma$ has a transient of at
most~$3n^2$ in the following way:
Take a node~$h$ of~$\Gamma$ and set $B=A^{\otimes \ell(\Gamma)}$.  As~$G(B)$ contains a
self-loop with maximal ratio at node~$h$, the sequences~$B_{i,h}^{\otimes k}$
and~$B_{h,j}^{\otimes k}$ have transients at most~$n-1$, which implies that the
sequences~$A_{i,h}^{\otimes k}$ and~$A_{h,j}^{\otimes k}$ have transients at
most $(n-1)\cdot \ell(\Gamma)$.
Both have period~$\ell(\Gamma)$ and ratio~$A(\Gamma)/\ell(\Gamma)$.
It is not hard to show that, in this case, their {\em max-plus convolution}
\[
\bigoplus_{k_1+k_2=k} A_{i,h}^{\otimes k_1} \otimes A_{h,j}^{\otimes k_2}
\]
has the same period and ratio, and a transient of at most $2\cdot(n-1)\cdot
\ell(\Gamma) +
\ell(\Gamma)\leq 2n^2-n$.
This convolution is equal to the sequence of maximum weights of walks of
length~$k$ from~$i$ to~$j$ that visit node~$h$.
Because the ratios of these sequences, for all~$h$ in~$\Gamma$, are equal
to~$A(\Gamma)/\ell(\Gamma)$ and their transients are at most~$2n^2-n$, the sequence of
maxima, formed over all~$h$ in~$\Gamma$, has the same ratio and a transient of at
most~$2n^2-n$.
This argument, which is essentially identical to the one given by Nachtigall,
yields a bound of~$2n^2-n$, improving the bound of~$3n^2$ in
Theorem~\ref{thm:nachtigall}.

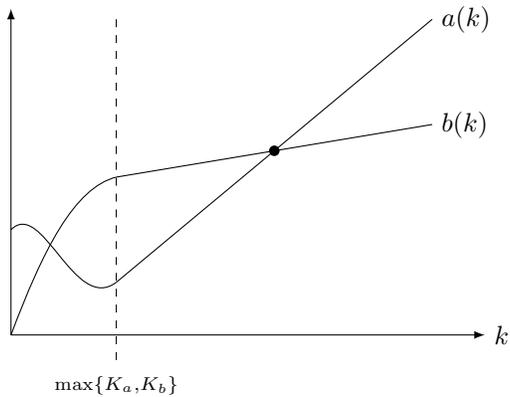
\begin{figure}
\centering
\begin{tikzpicture}[scale=1.4 , >=latex]
\draw[->] (0,0) -- (4.5,0) node [right] {$k$};
\draw[->] (0,0) -- (0,3.1) node [above] {};

\draw[dashed] (1,3) -- (1,-0.3) node [below] {$\scriptstyle \max\{K_a,K_b\}$};

\draw (0,1.0) .. controls (0.3,1.3) and (0.6,0.2) .. (1,0.5);
\draw (1,0.5) -- (4,3) node [right] {$a(k)$};

\draw (0,0.0) .. controls (0.3,0.8) and (0.6,1.4) .. (1,1.5);
\draw (1,1.5) -- (4,2.0) node [right] {$b(k)$};

\coordinate (t) at (intersection of 1,0.5 -- 4,3 and 1,1.5 -- 4,2.0);
\fill (t) circle (0.05) node [below right] {};

\end{tikzpicture}
\caption{Two eventually periodic sequences with differing ratios and respective
transients $K_a$ and $K_b$}
\label{fig:lines}
\end{figure}

We would like to point out that, if the ratios differ, the transient of a maximum of
eventually periodic sequences need not be bounded by the maximum of the
sequences' transients (see Figure~\ref{fig:lines}).
It is possible that the maximum is not even eventually periodic:
If~$a(k)$ and~$b(k)$ are two eventually periodic scalar sequences such that
$a(k)$'s ratio is strictly larger than that of~$b(k)$, then
the maximum $c(k) = a(k) \oplus b(k)$ is eventually periodic if and only if,
for all~$k$ large enough, $a(k)=\bzero$ implies~$b(k)=\bzero$.
This condition is not necessary for eventual periodicity if the two ratios are equal.

Thus, Theorem~\ref{thm:nachtigall} by itself, even if~$A$ is irreducible, i.e.,
$A^{\otimes k}$ is eventually periodic, does not give an upper bound on the
transient of~$A^{\otimes k}$.
However, inspection of its proof does:
By construction, if the $(i,j)$th entry of matrix~$A_r(k)$ is finite,
then 
there exists a walk of length~$k$ in~$G(A)$ from~$i$ to~$j$.
Because~$G(A)$ is strongly connected, by the finiteness of its index of convergence (see Section~\ref{sec:boolean}),
if~$k$ is large enough, there exists a walk of length~$k$ in~$G(A)$ from~$i$
to~$j$ that visits an arbitrary prescribed node~$h$.
This shows that, if~$k$ is large enough, then also the $(i,j)$th entry of~$A_1(k)$ is
finite.\footnote{
We do not need to quantify the threshold for~$k$ because eventual periodicity of~$A_1(k)$
and~$A_r(k)$ shows that the implication is true from $2n^2-n$ on, i.e., for
all $k\geq 2n^2-n$, as soon as we know that it is true from some (unknown) threshold on.
This suffices for our purposes.
}
Hence all sequences of the form~$A_1(k) \oplus A_r(k)$ are eventually periodic.
Their ratios are all equal to~$\lambda(A)$ since the first cycle to be removed
is a critical one. 
Hence the transient of the sequence $A^{\otimes k} = A_1(k)\oplus
A_2(k)\oplus \cdots \oplus A_n(k)$ is bounded by the maximum transient of
the $A_1(k)\oplus A_r(k)$'s.

By elementary calculations, one can show the following lemma which provides a bound on the transient
of the maximum of two eventually periodic sequences if the maximum itself is
eventually periodic.
It enables us to bound the transient of $A_1(k)\oplus A_r(k)$.

\begin{lemma}\label{lem:sequences}
Let~$a(k)$ and~$b(k)$ be two eventually periodic scalar sequences with
respective ratios $\varrho_a=0>\varrho_b$, transients at most $K$, and periods at
most~$p$ such that
$a(k)=\bzero$ implies $b(k)=\bzero$ for all $k$ large enough.
Then the maximum $a(k)\oplus b(k)$ is eventually periodic with transient at
most $K +p-1+ \Delta / (\varrho_a - \varrho_b)$ where $\Delta$ is the maximum
value of the form $b(k) - a(l)$ with $k,l\in \{ K, \dots , K+p-1\}$ and
$a(l)\neq \bzero$.
\end{lemma}

This lemma enables us to deduce the following transience bound from the proof
of the Nachtigall decomposition.
Without loss of generality, we can assume that $\lambda(A)=0$.
When applying the lemma to~$A_1(k)$ and~$A_r(k)$, the number~$K$ is at most
$2n^2-n$ and
$p$ is at
most~$n$, which implies that $\Delta$ is at most $2n^2\lVert
A\rVert$.

\begin{theorem}\label{thm:easy}
Let~$A$ be an irreducible $n\times n$ max-plus matrix.
Then the transient of the sequence of powers~$A^{\otimes k}$ is at most
\[
2n^2 + \frac{2n^2\lVert A\rVert}{\lambda(A) - \lambda_2(A)} 
\enspace.
\]
\end{theorem}

\section{Bound by Hartmann and Arguelles}\label{sec:ha}

Hartmann and Arguelles~\cite{HA99} gave the first general transience bound for
arbitrary irreducible max-plus matrices.
Their proof is purely graph-theoretic.

When analyzing their proof, one can extract a global proof strategy, variants of which are also found in later proofs of transience
bounds~\cite{SyK03,CBFN12}.
It has been described explicitly by Charron-Bost et al.~\cite[Section 3]{CBFN12}.
In order to prove that some number~$B$ is an upper bound on the transient of the
sequence~$A^{\otimes k}$ for an
irreducible
matrix~$A$, do the following:
\begin{enumerate}
\item Show that one can assume~$\lambda(A)=0$, i.e., the
sequence~$A^{\otimes k}$ is eventually periodic with
ratio~$0$.
\item Fix two nodes~$i$ and~$j$, and a congruence class~$[k_0]$
modulo some period~$p$ of the sequence~$A^{\otimes k}$. 
\item 
The assumption $\lambda=0$ guarantees that the maximum
$\max_{k\in M} A_{i,j}^{\otimes k}$ formed over an arbitrary nonempty
set~$M$ of nonnegative integers exists. 
We choose the set~$M$ to consist of those elements of class~$[k_0]$ that are greater or equal to~$B$.
Since the maximum exists, there exists a walk~$W$ from~$i$ to~$j$ with length
in~$M$ that attains it.
If~$B$ is indeed an upper bound on the transient,
the values $A^{\otimes k}_{i,j}$ with~$k\in M$ will all be equal.
\item Show that, whenever the length of~$W$ is greater or equal to some ``critical
bound'' $B_c\leq B$,
then it necessarily shares a node with a critical cycle~$\Gamma$.
\item 
Show that one can reduce walk~$W$ by removing subcycles such that it is
possible to attain all lengths
in~$M$ greater or equal to some ``pumping bound'' $B_p\leq B$ 
by adding
critical cycles.
The assumption $\lambda=0$ implies that all subcycles have weight at most~$0$
and critical cycles have weight
equal to~$0$.
Thus the weights of walks obtained in this way cannot be lower than that
of~$W$;
hence they are equal to that of~$W$.
\item We then have shown, because the choice of~$[k_0]$ was arbitrary, that the
transient of~$A^{\otimes k}_{i,j}$ is at
most $B\geq \max\{B_c,B_p\}$.
\end{enumerate}

Hartmann and Arguelles used $p=\gamma_c(A)$ in step~(2).
For step~(5), they described a walk reduction based
on the following basic application of the pigeonhole principle:

\begin{lemma}\label{lem:nt}
Let~$d$ be a positive integer.
Every collection of at least~$d$ integers has a nonempty subcollection whose
sum is divisible by~$d$.
\end{lemma}

They used this lemma to reduce walk~$W$ in step~(5).
After their reduction their walk could be disconnected, but they showed that
adding a copy of (critical) cycle~$\Gamma$ reestablishes connection
\cite[Theorem 4]{HA99}:

\begin{lemma}[Hartmann and Arguelles, 1999]\label{lem:ha:reduction}
Let~$W$ be a walk that shares a node with some cycle~$\Gamma$ and let $k$ be an
integer such that $k \equiv
\ell(W) \pmod{\ell(\Gamma)}$ and $k \geq n^2$ where~$n$ denotes the number of
nodes in the
graph.
Then there exists a walk~$\tilde{W}$ obtained from~$W$ by removing cycles and possibly adding
copies of\/~$\Gamma$ such that 
$\ell(\tilde{W})=k$.
\end{lemma}

To pump the walk length after the walk reduction,
they used
a result by Brauer~\cite{Bra42} on the Frobenius problem
to combine critical cycles to
attain a multiple of~$\gamma_c(A)$.
The use of Brauer's theorem introduces a term that is necessarily quadratic
in~$n$ to the
transience bound.
We want to note at this point that this use of Brauer's theorem can be avoided
by considering a period in step~(2) different from the critical subgraph's
cyclicity  because of Lemma~\ref{lem:transient:indep}.\footnote{
Hartmann and Arguelles actually prove
Lemma~\ref{lem:transient:indep} later in the paper
\cite[Lemma~11]{HA99}, but do not use it in the proof of their transience
bound.
}

The same strategy as described above can be adapted to show
transience bounds for systems $A^{\otimes k}\otimes v$.
In the case that all entries of~$v$ are finite, it is possible to show a sharper bound because
the walks under consideration do not have both the start and the end node
fixed, but only the start node.
This allows to circumvent the necessity of showing the existence of walks of
prescribed length between two fixed nodes (see
Section~\ref{sec:boolean}).

\begin{theorem}[Hartmann and Arguelles, 1999]\label{thm:ha}
Let~$A$ be an irreducible $n\times n$ max-plus matrix.
Then the transient of the sequence of powers~$A^{\otimes k}$ is at most
\[
\max \left\{ 2n^2 \ ,\ \frac{2n^2\lVert A \rVert}{\lambda(A) - \lambda_0(A)}  \right\}
\enspace.
\]
If, additionally, $v$ is a column
vector of size~$n$ with only finite entries,
then the transient of the system~$A^{\otimes k}\otimes v$ is at most
\[
\max \left\{ 2n^2 \ ,\ \frac{\lVert v \rVert + n\lVert A \rVert}{\lambda(A) - \lambda_0(A)}  \right\}
\enspace.
\]
\end{theorem}

Hartmann and Arguelles also proved a form of asymptotic tightness of their
transience bound for matrices.
They gave, for every~$n$ of the form $n=3m-1$ and all positive
reals~$\lambda$ and~$\lambda_0$ with~$\lambda>\lambda_0$, an irreducible $n\times n$
max-plus matrix~$A$ with~$\lambda(A)=\lambda$ and~$\lambda_0(A)=\lambda_0$ (see
\cite[Figure 1]{HA99}).
Their example has the property that $\lambda_0(A)=\lambda_{nc}(A)=\lambda_2(A)$
and $\lVert A\rVert = \lambda$.
They showed by explicit calculation that $A$'s transient is at
least $3 + m(m-2)\lambda/(\lambda - \lambda_0)$.

We can generalize their example to arbitrary~$n$ by inserting additional nodes that
do not change the transient.
This then shows that, even if one can prescribe all the other parameters in the
matrix bound of Theorem~\ref{thm:ha}, it is asymptotically tight when~$n$ tends
to infinity:

\begin{theorem}\label{thm:tight}
Let~$D_n$ and~$M_n$ be two sequences of positive real numbers such that $D_n
\leq M_n$.
Then there exists a sequence of irreducible $n\times n$ max-plus matrices~$A_n$
such that $\lambda(A_n) - \lambda_2(A_n) = D_n$, $\lVert A_n\rVert = M_n$, and
the transient of the sequence of matrix powers~$A_n^{\otimes k}$ is
\[
\Omega \left( \frac{  n^2  \lVert A_n\rVert}{\lambda(A_n) -
\lambda_2(A_n)} \right)
\enspace.
\]
\end{theorem}

Because $\lambda_2=\lambda_{nc}=\lambda_0$ in Hartmann and
Arguelles' example, Theorem~\ref{thm:tight} also holds with either
$\lambda_{nc}$ or
$\lambda_0$ replacing $\lambda_2$.

\section{A Bound for Primitive Matrices}\label{sec:bg}

A certain class of graph-theoretic arguments
 has been developed for the case
that the matrix is {\em primitive}, i.e., if its critical subgraph has a cyclicity equal
to~$1$.\footnote{This definition is consistent with the definition of
primitivity for Boolean
matrices (Section~\ref{sec:boolean}) because all cycles are critical in the
Boolean case.}
This class of arguments was used by both
Akian et al.~\cite[Remark 7.14]{AGW05}
and
Bouillard and Gaujal~\cite{BG01}.
To explicitly state a bound emerging from these arguments, we present the bound of Bouillard and Gaujal in
this section.
For ease of notation, we give it for the case~$\lambda(A)=0$.

\begin{theorem}[Bouillard and Gaujal, 2001]\label{thm:bg}
Let~$A$ be a primitive irreducible $n\times n$ max-plus matrix with
$\lambda(A)=0$.
Then the transient of the sequence of powers~$A^{\otimes k}$ is at most
\[
\max \left\{ 2n-2 + H + (n_c - 2H)\cdot\hat{g} \ ,\ \frac{\max_{i,j}
\lvert W_{i,j}^{(nc)} - W_{i,j}^{(c)} \rvert}{ -
\lambda_{nc}(A)} + (n-n_c)  \right\}
\]
where~$n_c$ is the number of critical nodes, $H$ is the number of critical
components, $\hat{g}$ is the maximum girth of components of the critical
subgraph, $W_{i,j}^{(nc)}$ is the maximum weight of walks from~$i$ to~$j$ not
visiting a critical node, and $W_{i,j}^{(c)}$ is the maximum weight of walks
from~$i$ to~$j$ that do visit a critical node.
\end{theorem}

Bouillard and Gaujal explained how to extend their result to the case of non-primitive 
matrices:
If~$A$'s critical subgraph has cyclicity~$\gamma_c$, then~$A^{\otimes\gamma_c}$ is
primitive. 
It is not necessarily irreducible, but it is guaranteed to be {\em completely
reducible}, i.e., permutation similar to a blockwise diagonal matrix
whose diagonal blocks are irreducible.
Also, every irreducible block contains at least one critical cycle, i.e., their
eigenvalues are equal, which implies that the sequence of powers is eventually
periodic.
If~$K$ is the transient of the sequence~$A^{\otimes k\gamma_c}$, then the
transient of~$A^{\otimes k}$ is at most~$K\gamma_c$.

Unfortunately, the cyclicity~$\gamma_c$ can be exponential in the size~$n$ of the
matrix.  This was shown by Malka et al.~\cite[Theorem 4]{MMZ93} who constructed matrices
whose critical subgraphs are disjoint unions of cycles of prime lengths. 
Using the Prime Number Theorem, one sees that it is possible to construct a
critical subgraph with cyclicity~$\gamma_c = e^{\Omega(\sqrt{n})}$.
Malka et al.\ improved this observation by showing that even
the {\em minimal\/} period can be in the same order:

\begin{theorem}[Malka et al., 1993]
There exists a sequence of irreducible $n\times n$ max-plus matrices~$A_n$ such
that the minimal
period of the sequence of matrix powers~$A_n^{\otimes k}$ is\/ $\exp\!\big( \Omega(\sqrt{n}) \big)$.
\end{theorem}

\section{When All Entries Are Finite}\label{sec:syk}

Soto y Koelemeijer~\cite[Theorem 3.5.12]{SyK03} established a 
transience bound
in the case that all matrix entries are finite, i.e., the corresponding
digraph is the complete graph.
His approach is similar to that of Hartmann and Arguelles, but the assumption
of existence of all edges in the corresponding digraph allows to construct
shorter walks.
Utilizing this fact, he arrived at a bound that can be lower than that of Hartmann and
Arguelles
(first part of Theorem~\ref{thm:ha}).

\begin{theorem}[Soto y Koelemeijer, 2003]
Let~$A$ be an $n\times n$ max-plus matrix with only finite entries.
Then the transient of the sequence of powers~$A^{\otimes k}$ is at most
\[
\max\left\{2n^2 \ ,\ \left\lceil\frac{ 2 \lVert A \rVert}{\lambda(A) -
\lambda_2(A)}\right\rceil + n-1\right\}
\enspace.
\]
\end{theorem}

\section{Inclusion of Cyclicity and Girth}\label{sec:cbfn}

Charron-Bost et al.~\cite{CBFN12} gave two transience bounds for systems; one
that involves the maximum girth of connected components of the critical
subgraph, and one that involves the maximum cyclicity.

For both, they used the general proof
strategy that we detailed in Section~\ref{sec:ha}.
With respect to this strategy, they
introduced two different methods of deleting and adding cycles in step~(5).
Their use is enabled by considering a period~$p$ equal to the least common
multiple of all critical cycle lengths in step~(2).

The two methods are called the {\em repetitive\/} and the {\em explorative\/}
method.
Denote by~$h$ the a critical node of walk~$W$ in step~(5).
Depending on the strongly connected component of the critical subgraph that~$h$
lies in, they choose a positive integer~$d$ as a parameter for the walk
reduction~$\Red_{d,h}$, whose definition we present later.
It has the following properties:

\begin{lemma}[Charron-Bost et al., 2012]\label{lem:red:upper:bound}
Let~$W$ be a walk containing node~$h$ and let~$d$ be a positive integer.
Then there exists a walk $\hat{W}=\Red_{d,h}(W)$ obtained from~$W$ by removing
subcycles such that
(i) $\tilde{W}$ still contains node~$h$, (ii) $\ell(\tilde{W}) = \ell(W)
\pmod{d}$, and (iii)
$\ell(\hat{W}) \leq 2\cdot d\cdot (n-1) + d- 1$, where~$n$ denotes the number
of nodes in the digraph.
\end{lemma}

In the repetitive method, they choose some critical cycle~$\Gamma$ which node~$h$ is
part of and choose~$d=\ell(C)$.
For pumping the reduced walk, they add copies of~$\Gamma$.
Since~$d$ divides~$p$, all lengths in the congruence class $[\ell(W)]$
modulo~$p$ can be reached that are greater or equal to~$\ell(\hat{W})$.

In the explorative method, they choose~$d$ equal to the cyclicity of $h$'s
strongly connected component in the critical subgraph.
For pumping in the explorative method, they add a closed walks in
the component starting at~$h$ using the notion of index of convergence (see
Section~\ref{sec:boolean}).
Again,~$d$ divides~$p$.
Hence all lengths in the congruence class $[\ell(W)]$
modulo~$p$ can be reached that are greater or equal to~$\ell(\hat{W}) + \ind$,
where~$\ind$ denotes the component's index of convergence.

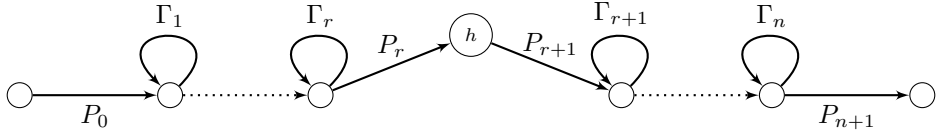
\begin{figure}
\centering
\begin{tikzpicture}[>=latex',scale=1.0]
	\node[shape=circle,draw] (i) at (-6,0) {};
	\node[shape=circle,draw] (k) at (0,.8) {$\scriptstyle h$};
	\node[shape=circle,draw] (j) at (6,0) {};
	\node[shape=circle,draw] (n1) at (4,0) {};
	\node[shape=circle,draw] (n2) at (2,0) {};
	\node[shape=circle,draw] (n3) at (-2,0) {};
	\node[shape=circle,draw] (n4) at (-4,0) {};
	\draw[thick,->] (i) -- node[below] {$P_0$} (n4);
	\draw[thick,->,dotted] (n4) -- node[below] {} (n3);
	\draw[thick,->] (n3) -- node[above] {$P_r$} (k);
	\draw[thick,->] (k) -- node[above] {$P_{r+1}$} (n2);
	\draw[thick,->,dotted] (n2) -- node[below] {} (n1);
	\draw[thick,->] (n1) -- node[below] {$P_{n+1}$} (j);
	\draw[thick,->] (n4) .. controls +(1,1) and +(-1,1) .. node[above]
{$\Gamma_1$} (n4);
	\draw[thick,->] (n3) .. controls +(1,1) and +(-1,1) .. node[above]
{$\Gamma_r$} (n3);
	\draw[thick,->] (n2) .. controls +(1,1) and +(-1,1) .. node[above]
{$\Gamma_{r+1}$} (n2);
	\draw[thick,->] (n1) .. controls +(1,1) and +(-1,1) .. node[above]
{$\Gamma_n$} (n1);
\end{tikzpicture}
\caption{Structure of the reduced walk $\hat{W} = \Red_{d,h}(W)$}
\label{fig:lem:red:upper:bound}
\end{figure}

The walk reduction $\Red_{d,h}(W)$ is defined as follows:
Starting at~$W$, one repeatedly removes nonempty collections of subcycles such that (a)
their combined length is a multiple of~$d$ and (b) after their removal, $h$ is
still a node of the resulting walk.
If there are more than one such collections, choose any.
Eventually, all collections of subcycles that satisfy (a) and (b) will be
empty.
At this point, the walk reduction stops and returns the last walk.
The resulting walk then has a structure as depicted in
Figure~\ref{fig:lem:red:upper:bound}:
It is a sequence of paths~$P_t$ connecting the start node to the end node via
some cycles~$\Gamma_t$ and the node~$h$.
By using Lemma~\ref{lem:nt}, one sees that there can be at most~$d-1$
cycles~$\Gamma_t$ because otherwise they could be removed as they would satisfy (a)
and (b).
It follows that there are at most~$d+1$ paths~$P_t$.
Since the paths have length at most~$n-1$ and the cycles at most~$n$,
the upper bound of Lemma~\ref{lem:red:upper:bound} on the reduced walk length follows.

\begin{theorem}[Charron-Bost et al., 2012]
Let~$A$ be an irreducible $n\times n$ max-plus matrix and let~$v$ be a column
vector of size~$n$ with only finite entries.
Then the transient of the system~$A^{\otimes k}\otimes v$ is less or equal to both
\[
\max \left\{ 2\hat{g}\cdot (n-1) + \hat{g} - 1 \ ,\ \frac{\lVert v \rVert + (n-1)\lVert A \rVert}{\lambda(A)
- \lambda_{nc}(A)}  \right\}
\]
and
\[
\max \left\{ \hat{\ind} +  2\hat{\gamma}\cdot (n-1) + \hat{\gamma} - 1 \ ,\ \frac{\lVert v \rVert + (n-1)\lVert A \rVert}{\lambda(A)
- \lambda_{nc}(A)}  \right\}
\enspace,
\]
where $\hat{g}$, $\hat{\gamma}$, and\/ $\hat{\ind}$ denote the greatest girth,
cyclicity, and index of critical components of~$G(A)$, respectively.
\end{theorem}

\section{From Systems to Matrices}\label{sec:sys:mat}

Charron-Bost et al.~\cite{CBFN12} also showed how to transform bounds for system transients
into bounds for matrix transients.
They used the following idea:
The transient of the sequence of powers of a matrix~$A$ is equal to the maximum
transient of systems $A^{\otimes k}\otimes v$ where~$v$ is one of the max-plus
unit vectors.
However all transience bounds for systems assume the vector~$v$ to have only
finite entries, which is not the case for the max-plus unit vectors.
So they considered ``truncated'' unit vectors that have their infinite entries
replaced by $-\mu$ where~$\mu$ is an appropriately chosen real number.
They used this approach, together with graph-theoretical arguments, to show
the following theorem.

\begin{theorem}[Charron-Bost et al., 2012]
Let~$A$ be an irreducible $n\times n$ max-plus matrix and let~$B$ be a
nonnegative integer.
Set
\[
\tilde{B} = 2n - 3 + \hat{\ind} + \ind\!\big(G(A)\big) + \hat{\gamma} 
\]
where $\ind\!\big(G(A)\big)$ denotes the index of convergence of~$G(A)$, and $\hat{\ind}$ and $\hat{\gamma}$ denote the greatest index and cyclicity
of strongly connected components of the critical subgraph.

If~$B$ is an upper bound on all transients of systems~$A^{\otimes k}\otimes v$
with $\lVert v\rVert \leq \tilde{B}\cdot \lVert A\rVert$, then the transient of
the sequence of matrix powers~$A^{\otimes k}$ is at most~$\max\{B,\tilde{B}\}$.
\end{theorem}

\section{Conclusion}\label{sec:conclusion}

We have presented various transience bounds and some of their proofs for both
max-plus systems
and matrices.
Most of the proofs were heavily graph-theoretic, with the exception of 
Theorem~\ref{thm:easy}, which is more algebraic and founded
on the concept of convolution of sequences.
Except for the bound of Theorem~\ref{thm:easy},
which can be seen to be strictly greater than the others, there is no general
ordering between pairs of bounds.
This is due to the fact that all of them consider either a different set of
parameters or a special case.
A ``good'' choice of parameters is not an obvious to make.
At one extreme, one could declare the transient itself as a parameter, which
would lead to a trivial bound.
At the other extreme, restricting oneself to only consider the matrix
size~$n$ as a parameter is not tractable either. 
In fact, the tightness result of
Theorem~\ref{thm:tight} shows that no upper bound only in terms of~$n$ exists.
We think that the choice of parameters can only depend on the envisioned
application of the respective transience bound.
Because the problem of computing the exact transient is computationally
feasible, parameters should be ones that can be controlled
during the system design phase.

\section*{Acknowledgments}

The authors would like to thank Marianne Akian and Anne Bouillard.

\bibliography{paper}
\bibliographystyle{amsplain}

\end{document}